\documentclass[11pt,leqno]{article}
\usepackage{amsmath,amssymb,amsthm}
\usepackage[mathscr]{eucal}
\DeclareMathOperator{\dist}{dist}
\DeclareMathOperator{\Span}{Span}

\newtheorem{thm}{Theorem}[section]

\newtheorem{lem}{Lemma}[section]
\newtheorem{prop}{Proposition}[section]
\newtheorem{cor}{Corollary}[section]

\theoremstyle{definition}
\newtheorem{defn}{Definition}[section]
\newtheorem{conj}{Conjecture}[section]

\newtheorem{rem}{Remark}[section]

\newcommand{\Rd}{\mathbb{R}^d}
\newcommand{\Rt}{\mathbb{R}^2}

\newcommand{\R}{\mathbb{R}}
\newcommand{\N}{\mathbb{N}}
\newcommand{\NP}{\mathbb{N}}

\newcommand{\eps}{\varepsilon}

\numberwithin{equation}{section}

\begin{document}
\newcommand{\bD}{\mathrm{I\! D\!}}

\newcommand{\al}{\alpha}
\newcommand{\ga}{\gamma}
\newcommand{\om}{\omega}
\newcommand{\G}{\Gamma}
\def\rayo{\leftarrow}
\def\tr{\triangle}
\newtheorem{cl}{Claim}
\newcommand{\bR}{\mathbf{R}}

\theoremstyle{definition}
\newtheorem{fact}{Fact} \renewcommand{\thefact}{}

\begin{titlepage}
\title{\bf  Eigenvalue  gaps for the Cauchy process and a Poincar\'e inequality}
\author{Rodrigo Ba\~nuelos\thanks{Supported in part by NSF Grant
\# 9700585-DMS}\\Department of  Mathematics
\\Purdue University\\West Lafayette, IN
47906\\banuelos@math.purdue.edu\and Tadeusz  Kulczycki\thanks{Supported by KBN 
grant 2 P03A 041 22
and RTN Harmonic Analysis and  Related Problems, contract 
HPRN-CT-2001-00273-HARP}
\\Institute of Mathematics\\ Wroc{\l}aw
University of Technology\\ 50-370 Wroc{\l}aw, Poland\\
tkulczyc@im.pwr.wroc.pl}
\date{}
\maketitle

\begin{abstract}
{\it   A connection  between
the semigroup of the Cauchy process killed upon exiting a domain $D$ and a mixed boundary
value problem for the Laplacian in one dimension higher known as 
{\bf the mixed Steklov
problem}, was established in \cite{BK}.  From this, a variational characterization for
the eigenvalues $\lambda_n$, $n\geq 1$,  of the Cauchy process in $D$ was obtained. 
 In this paper
 we obtain a variational
characterization of the difference between  $\lambda_n$ and $\lambda_1$.  We study
bounded convex domains which are symmetric with respect to one of the coordinate axis  and
obtain lower bound estimates for $\lambda_* - \lambda_1$ where $\lambda_*$ is the
eigenvalue corresponding to the ``first'' antisymmetric eigenfunction for $D$. 
The proof is based on a 
variational characterization of $\lambda_* - \lambda_1$ and on a weighted 
Poincar\'e--type inequality.  The Poincar\'e  inequality is  valid for all  $\alpha$
symmetric stable processes, $0<\alpha\leq 2$,  and  any other process obtained from
Brownian motion by subordination.  We also  prove upper bound estimates for the spectral
gap $\lambda_2-\lambda_1$ in bounded  convex domains.}

\smallskip

\centerline{\bf  Contents}
\begin{itemize}
\item[{\S1.}] {\sl Introduction}
\item[\S2.] {\sl  Variational formulas}
\item[\S3.] {\sl Weighted Poincar\'e inequalities}
\item[\S4.] {\sl Proof of Theorems \ref{main} and \ref{upperestimate}}
\item[\S5.] { \sl Concluding Remarks}
\end{itemize}
\end{abstract}
\end{titlepage}

\section{Introduction} 
\label{sec:introduction}

  The spectral gap estimates for eigenvalues of the Laplacian with Dirichlet
boundary conditions, henceforth referred to as the {\it Dirichlet Laplacian,}  have attracted
considerable  attention for many  years
\cite{AB1},
\cite{AB2}, \cite{vanden}, \cite{L}, \cite{Sm}, \cite{SWY}, \cite{YZ} . The Dirichlet
Laplacian is the infinitesimal generator of the semigroup of Brownian
motion killed upon leaving a domain. Therefore questions concerning
eigenvalues of this operator can be studied both by analytic and
probabilistic methods.  The question of  precise  lower bounds for  the spectral gap for
the Dirichlet Laplacian (the difference between the first two eigenvalues)
was raised  by M. van den Berg \cite{vanden} (see also Yau \cite{Yau}, problem \#44) and was
motivated by problems in mathematical physics related to the behavior of free
Boson gases. 
The conjecture, which remains open,  asserts that for any convex bounded domain $D$ of
diameter
$d_D$, the spectral gap is bounded below by $3\pi^2/d_D^2$. (See \cite{BKr}, \cite{BM},
\cite{Da} where some special cases of the conjecture are proved and \cite{DH}, \cite{You} for
more genral ``partition  function" inequalities.)   The spectral gap has also been studied for the
Laplacian with Neumann boundary conditions and for Schr\"odinger operators,
\cite{PW}, \cite{Sm},
\cite{AB1}, \cite{SWY}. 
 From the probabilistic point of view, the
spectral gap for the Dirichlet Laplacian determines the  rate to equilibrium for  the Brownian
motion conditioned to remain forever in
$D$, the Doob \emph{h--process} corresponding to the ground state eigenfunction. 

The natural question arises as to 
 whether these results can
be extended to other  non-local, pseudo-differential operators. The class of such
operators which are most closely related   to the Laplacian $\Delta$ from the point of view of
Brownian motion are $-(-\Delta)^{\alpha/2}$, $\alpha \in
(0,2)$.  These are   the infinitesimal
generators of the symmetric $\alpha$-stable processes.  
 These processes do not have continuous paths which is related
to non-locality of $-(-\Delta)^{\alpha/2}$. As in the case of Brownian
motion, we can consider the semigroup of these 
processes killed upon exiting domains and we can consider the eigenvalues 
 of such semigroup. Here again, the spectral gap
determines the asymptotic exponential rate of convergence to
equilibrium for the process conditioned to remain forever in the domain. 
Instead of speaking of
the eigenvalue gap for the operator  $-(-\Delta)^{\alpha/2}$ 
we will very often refer
to it as the eigenvalue gap for the corresponding process. 

The purpose of this paper is to obtain  eigenvalue gap estimates for the Cauchy
process, the symmetric $\alpha$-stable process for $\alpha
=1$. This is done using the connection (established in \cite{BK})
between the eigenvalue problem for the Cauchy process and 
 the {\it mixed Steklov problem.} Both, the  methods and the  results,
are new.  The results raise natural questions
concerning  spectral gaps for other
symmetric
$\alpha$-stable processes and for more general  Markov processes.  We  believe that
as with the results in \cite{BK} which have motivated subsequent work by others, 
see \cite {Bl2}, \cite{DeMe} \cite{CS3}, the current results will also be of interest.   
Let $X_t$ be a symmetric $\alpha$-stable process in $\Rd$, $\alpha \in
(0,2]$. This is a process with independent and stationary increments
and characteristic function 
$E^0 e^{i \xi X_t} = e^{-t |\xi|^{\alpha}}$, $\xi \in \Rd$, $ t > 0$. 
$E_x$, $P_x$ denote the expectation and probability of this process
starting at $x$, respectively.  By $p^{(\alpha)}(t,x,y) = p_t^{(\alpha)}(x-y)$ we
will denote the transition density of this process. That is, 
$$P_{x}(X_t \in
B)= \int_{B}p^{(\alpha)}(t,x,y), dy.
$$
When $\alpha = 2$ the process $X_t$ is just the Brownian motion in $\Rd$ but running at twice the speed.
That is, if $\alpha = 2$ then
\begin{equation}
\label{gaussian}
p^{(2)}(t,x,y) = \frac{1}{(4 \pi t)^{d/2}} e^{\frac{-|x-y|^2}{4t}},\quad t > 0, \, \, x,y \in \Rd.
\end{equation}
When $\alpha = 1$, the process $X_t$ is the Cauchy process in $\Rd$ whose transition densities are given
by
\begin{equation}
\label{transition}
p^{(1)}(t,x,y) = \frac{c_d \, t}{(t^2 + |x - y|^2)^{(d + 1)/2}}, \quad t > 0, \, \, x,y \in \Rd,
\end{equation}
where
$$
c_d = \Gamma((d+1)/2)/ \pi^{(d+1)/2}.
$$

Our main  interest in this paper are the eigenvalues of the semigroup of the process
$X_t$ killed upon  leaving a domain. Let $D \subset \Rd$ be a bounded connected domain
and $\tau_{D} =
\inf\{t
\ge 0: X_t
\notin D\}$ be the first exit time of $D$.   By $\{P_t^{D}\}_{t \ge 0}$ we denote the semigroup on
$L^2(D)$ of $X_t$ killed upon exiting $D$.  That is, 
 
$$
P_{t}^{D}f(x) = E_{x}(f(X_{t}), \tau_{D} > t), \quad x \in D, \, \, t > 0, \, \, f \in L^2(D).
$$
The semigroup has transition densities $p_D(t,x,y)$ satisfying 
$$
P_{t}^{D}f(x) = \int_{D} p_{D}(t,x,y) f(y) \, dy.
$$
The  kernel $p_{D}(t,x,y)$  is strictly positive
symmetric and
$$
p_{D}(t,x,y) \le p^{(\alpha)}(t,x,y) \le c_{\alpha,d} \, t^{-d/\alpha}, \quad x,y \in D, \, \, t>0.
$$
The fact that $D$ is bounded implies that for any $t > 0$ the operator $P_t^D$ maps $L^2(D)$ into $L^{\infty}(D)$.
From the
general theory of
 semigroups
\cite{Da1} it follows that there exists 
 an   orthonormal basis of eigenfunctions
$\{\varphi_n\}_{n =1}^{\infty}$ for
$L^2(D)$ and corresponding eigenvalues
$\{\lambda_n\}_{n = 1}^{\infty}$ satisfying
$$0<\lambda_1<\lambda_2\leq \lambda_3\leq \dots$$
  with $\lambda_n\to\infty$
as
$n\to\infty$. That is,  the pair $\{\varphi_n, \lambda_n\}$ satisfies
\begin{equation}
\label{cauchyproblem}
P_{t}^{D}\varphi_{n}(x) = e^{-\lambda_{n} t} \varphi_{n}(x), \quad x \in D, \,\,\, t > 0.
\end{equation}
The eigenfunctions $\varphi_n$ are continuous and bounded on $D$. In addition, 
$\lambda_1$ is simple and the corresponding eigenfunction
$\varphi_1$, often called the ground state eigenfunction, is strictly
positive on
$D$.  
 For more general properties of
  the semigroups $\{P_t^D\}_{t \ge 0}$, see
\cite{G1},
\cite{BG2}, \cite{CS1}.

It is well known (see \cite{B}, \cite{CS1}, \cite{CS2}, \cite{K})
that if $D$ is a bounded connected Lipschitz domain and $\alpha = 2$, or that if 
$D$ is a bounded connected domain for $0 < \alpha < 2$, then
$\{P_t^D\}_{t \ge 0}$ is intrinsically ultracontractive.  This implies, among many other things,
that  
$$
\lim_{t \to \infty} \frac{e^{\lambda_1 t} p_D(t,x,y)}{\varphi_1(x) \varphi_1(y)} = 1, 
$$
uniformly in both variables $x, y\in D$.  In addition, the rate of convergence is given by 
the spectral gap $\lambda_2-\lambda_1$. That is, 
for any $t \ge 1$ we have
\begin{equation}
\label{spectralgap}
e^{-(\lambda_2 - \lambda_1) t} \le 
\sup_{x,y \in D} \left|\frac{e^{\lambda_1 t}p_D(t,x,y)}{ \varphi_1(x) \varphi_1(y)} - 1
 \right| \le C(D,\alpha) e^{-(\lambda_2 - \lambda_1)t}.
\end{equation}
The proof of this  for $\alpha=2$ may be found in \cite{Sm}.  
The proof in our setting is exactly the same.

In the Brownian motion case the properties of
eigenfunctions  and eigenvalues have
been extensively
studied for many years, both analytically and probabilistically.
It is  well known
  that  geometric information
on $D$, such as convexity, symmetry, volume growth, smoothness of its
boundary, etc.,
provides   information not only on the ground state
eigenfunction $\varphi_{1}$
and the ground state  eigenvalue
     $\lambda_{1}$, but also on the  spectral gap
$\lambda_2-\lambda_1$,  and on the geometry of the nodal domains of
$\varphi_2$.

In the case of stable processes of index $0<\alpha<2$, very little is
known. (We refer the
reader to \cite{BK} where some of the known results are reviewed and for a discussion of the
many open questions.) Except for the one-dimensional case (\cite{BK}, \cite{CS3}) we are not
at present able to estimate from below  the spectral gap
$\lambda_2 -
\lambda_1$ or obtain much useful geometric  information on  the 
eigenfunction corresponding to $\lambda_2$.  In this paper we will instead study domains
with one  axis of symmetry and obtain estimates for $\lambda_* - \lambda_1$ where
$\lambda_*$ is the eigenvalue corresponding to the ``first'' antisymmetric eigenfunction
for $D$. In the Brownian motion case $\lambda_*=\lambda_2$
in many important cases (we will discuss this later in the
sequel). Therefore estimates on $\lambda_* -
\lambda_1$ are  very closely related to estimates on $\lambda_2 -
\lambda_1$. It is natural to conjecture that $\lambda_* - \lambda_1
= \lambda_2 - \lambda_1$ for the Cauchy process and for other symmetric $\alpha$-stable
processes in various symmetric domains but this remains open.

For each $x = (x_{1},x_{2},\ldots,x_{d})$ we put $\widehat{x} =
(-x_{1},x_{2},\ldots,x_{d})$.
For any domain $D\subset \Rd$, we set $D_{+} = \{x \in D: x_{1} >
0\}$  and $D_{-} = \{x \in D: x_{1} < 0\}$.
We say that $D$ {\it{ is symmetric relative 
to the $x_1$--axis}} if
$\widehat{x}\in D$ whenever $x\in D$. 
Recall that the inradius $r_D$ of $D$ is the radius of the largest ball
contained in $D$.

In Theorem 4.3 of \cite{BK}  we proved that 
if  $D \subset \Rd$ is  a connected, bounded Lipschitz domain which is symmetric
relative to the  $x_1$--axis, then there exists
an eigenfunction $\varphi_{*}$ for the Cauchy process
with corresponding eigenvalue $\lambda_{*}$
which   is antisymmetric relative to  the
$x_{1}$--axis ( $\varphi_{*}(x) = -\varphi_{*}(\widehat{x})$, $x \in D$) and
(up to a sign) $\varphi_{*}(x) > 0$ for $x \in D_{+}$ and $\varphi_{*}(x) < 0$
for $x \in D_{-}$. Moreover, if $\varphi$ is any  eigenfunction with
  eigenvalue $\lambda$ such that $\varphi$ is antisymmetric
relative to the  $x_{1}$--axis and $\varphi$ is different from 
$\varphi_{*}$ ($\varphi
\notin \Span\{\varphi_{*}\}$), then $\lambda_{*} < \lambda$. In other words,
$\varphi_{*}$ has the smallest eigenvalue amongst all eigenfunctions which are
antisymmetric relative to $x_{1}$--axis.

 The main result of this paper is the following theorem.

\begin{thm}
\label{main}
Let $D \subset \Rd$ be a bounded convex Lipschitz domain which is
symmetric relative to the $x_1$-axis and $\{P_t^D\}_{t \ge 0}$ be the
semigroup of the Cauchy process killed upon exiting $D$. Let $\lambda_*$ be
 the eigenvalue for $\{P_t^D\}_{t \ge 0}$ corresponding to the unique eigenfunction $\varphi_*$ which is
antisymmetric relative to the $x_1$-axis and strictly positive on $D_+$ and strictly negative on $D_-$.  Let
$L = \sup\{x_1: x=(x_1,\ldots,x_d) \in D\}$ and assume that the inradius $r_D$ of $D$ is equal
to
$1$. Then we have
\begin{equation}\label{lower}
 \min \left( \frac{C_d}{L^2}, C'_d \right)\leq \lambda_* - \lambda_1 
\end{equation}
where
$$
C_d =  \frac{\pi^2 (d+1)}{2 \pi \, d (d+2) + 4 (d+1)},
$$
and $C'_d = 4 C_d/\pi^2$.
\end{thm}

The eigenvalues $\lambda_n$ satisfies the scaling property
 $\lambda_n(kD) = \lambda_n(D)/k$, $k > 0$. This leads to the following easy conclusion.

\begin{cor}
Let $D\subset \Rd$ satisfy the same assumption as in Theorem \ref{main} except that
 now the inradius $r_D$ is arbitrary. Then we have
\begin{equation}
 \min \left(\frac{ C_d \, r_D}{L^2}, \frac{C'_d}{r_D} \right)\leq \lambda_* - \lambda_1 
\end{equation}
where $C_d$, $C'_d$ are the same as in Theorem \ref{main}.
In particular,  for a disk $D = B(0,r) \subset \Rt$, $r>0$ we have
$$
 \frac{1}{6 \, r}\leq \lambda_* - \lambda_1. 
$$ 
\end{cor}

In terms of an upper estimate for the gap,  we have the following 
\begin{thm}\label{upperestimate}
Let $D\subset \Rd$ be a bounded convex domain of inradius $r_D$ and let 
$\lambda_1$,
$\lambda_2$ eigenvalues for the semigroup of the Cauchy process killed upon exiting $D$. 
Then 
\begin{equation}\label{upper}
 \lambda_2 - \lambda_1 \leq
{\frac{\sqrt{\mu_2}-{\frac{1}{2}}\sqrt{\mu_1}}{r_D}}
\end{equation}
where $\mu_1$ and $\mu_2$ are, respectively, the first and second
 eigenvalues for the Dirichlet Laplacian for the unit ball, $B(0, 1)$, in $\Rd$.  In fact, $\mu_2=j_{{d/2}, 1}^2$ and 
 $\mu_1=j_{{d/2-1}, 1}^2$ where $j_{p, k}$ denotes the $k^{th}$ positive zero of the Bessel
function
$J_p(x)$. 
\end{thm}

The constants $C'_d$, $C'_d$ in Theorem \ref{main} are of course not optimal. An easy
calculation shows that $C_1
\approx 0.735$, ($C_1 > 7/10$),
$C'_1 \approx 0.297$ ($C'_1 > 1/4$),
$C_2 \approx 0.475$, ($C_2 > 4/10$), 
 $C'_2 \approx 0.192$, ($C'_2 > 1/6$), $C_3 \approx 0.358$, ($C_3 > 1/3$), $C'_3 \approx 0.145$, ($C'_3 >
1/7$).  In particular, for rectangles $R = [-L,L] \times [-1,1]$, where $L \ge 1$, we have
$$
 \min \left( \frac{2}{5 \, L^2}, \frac{1}{6} \right)<\lambda_* - \lambda_1 
$$

 As we shall see below, Theorem \ref{upperestimate}  holds in greater generality 
and it also raises 
interesting questions concerning sharp upper bounds; see Conjecture \ref{PPW} below.

 In the case of  Brownian motion under the same
assumptions on $D$ there is  also an  antisymmetric eigenfunction
$\varphi_*$. In fact, for Brownian motion $\varphi_*$ restricted to
$D_+$ is  the first eigenfunction for $D_+$ and hence $\lambda_*(D) =
 \lambda_1(D_+)$. This fact has been used by several authors
to study  the van den Berg conjecture  mentioned above (\cite{BM}, \cite{Da}).   For the
Cauchy process
$\lambda_1(D_+)
\ne
\lambda_*(D)$ 
(in fact $\lambda_1(D_+) < \lambda_*(D)$) 
and
$\varphi_*$ restricted to $D_+$ is not the  first
eigenfunction for $D_+$. Such effect is due to the discontinuity of the paths of the 
Cauchy process. 
This is the reason for
introducing the special eigenvalue $\lambda_*$ instead of studying
$\lambda_1(D_+)$ as in the case of Brownian motion.

In the case
 of Brownian motion for a bounded domain $D$ the
 Courant-Hilbert nodal domain theorem asserts that the second
 eigenfunction $\varphi_2$ has exactly 2 nodal domains. That is, $D$
 is divided into 2 connected subdomains $D_+$ and $D_-$ such that
 $\varphi_2 > 0$ on $D_+$ and $\varphi_2 < 0$ on $D_-$.
If in addition  $D \subset \Rt$ is convex,  the nodal
line $N =
 \{x \in D: \, \varphi_2(x) = 0\}$ touches the boundary at exactly 2
 points (\cite{M}, \cite{A}). Moreover, when $D \subset \Rt$ is convex 
 and double symmetric, that is, $D$ is symmetric relative to both
 coordinate axes,  there exists an  eigenfunction corresponding
 to $\lambda_2$ with nodal line lying on one of the coordinate axes
 (see L. E. Payne \cite{P1}). In other words $\varphi_2 = \varphi_*$
 or $\varphi_2$ is  an  antisymmetric eigenfunction defined
 analogously as $\varphi_*$ but with respect  to the $x_2$--axis. Therefore in
 the case of Brownian motion, when $D \subset \Rt$ is a convex double
 symmetric domain,
 estimates for $\lambda_* - \lambda_1$ gives estimates for $\lambda_2
 - \lambda_1$. 
Unfortunately,  in the case of the Cauchy
process we do not know anything about the location of 
  the nodal line for the second eigenfunctions even in the simplest possible planar regions
such as  a  disk or  a rectangle. Nevertheless, it seems that the
following conjecture should be true.

\begin{conj}
\label{hyp1}
Let $D \subset \Rt$ be a convex   domain which is
symmetric relative to both coordinate axis. Let $\lambda_n$,
$\varphi_n$ be the eigenvalues and eigenfunctions 
for the Cauchy process in $D$.  Then there exists an eigenfunction
 corresponding to $\lambda_2$ such that its nodal line lies on one of the coordinate axis. 
\end{conj}

 If this conjecture were  true then the estimates for $\lambda_* -
\lambda_1$ would give estimates for $\lambda_2 - \lambda_1$. 
We are not able to prove the Conjecture \ref{hyp1} partly because we do not know whether the
 Courant-Hilbert nodal domain theorem holds for the Cauchy process. 
 It may be possible to gain some information on this conjecture 
 by analyzing $\partial \varphi_1(x)/\partial x_i$  as in \cite{P1} but so far this remains open. 
In the simplest geometric situation of  $D = (-1,1)$ we know the ``shape" of the second
eigenfunction, that  $\lambda_2$ has
multiplicity $1$  and  that $\lambda_2 = \lambda_*$ , (\cite{BK}, Theorem 5.3). However,
even in this simple geometric setting the  situation is fairly nontrivial.

The paper is organized as follows. In \S2, we recall  the
connection between eigenvalues and eigenfunctions for the
Cauchy process and the Steklov problem, (\cite{BK}).   Using this we
derive a variational formulas for $\lambda_* - \lambda_1$ and $\lambda_n
- \lambda_1$. Such variational formulas are of independent
interest. Also, in \S2 we present some auxiliary  lemmas which allow us to replace the Steklov eigenfunction $u_1(x,t)$ in the variational formula by the simpler expression $e^{- \lambda_1 t} \varphi_1(x)$.

In \S3, we prove the weighted Poincar\'e--type inequality for the first eigenfunction $\varphi_1$. The
Poincar\'e inequality has been used in the Brownian motion case  in \cite{SWY} and \cite{Sm} 
 to estimate $\lambda_2
-\lambda_1$.  In that case the Poincar\'e  inequality  depends on the fact that 
for convex domains the first eigenfunction $\varphi_1$ is log--concave. 
For the  Cauchy process this remains unknown. (For some geometric properties related to
concavity for the eigenfunction in rectangles, see \cite{BKM}.)
 Nevertheless, by subordination we can show that
$\varphi_1$ is the limit of integrals of log-concave functions and this 
allows us to obtain the appropriate inequality. We will show this Poincar\'e inequality not only for the Cauchy
process but for all symmetric $\alpha$-stable processes $0 < \alpha < 2$.

In \S4 we prove  Theorems \ref{main} and \ref{upperestimate}. The lower bound (Theorem \ref{main})  will follow from the variational formula and the Poincar\'e inequality.  The upper bound is an easy observation
that follows from a deep result of Ashbaugh and Benguria, \cite{AB2},  and a recent result
of Chen and Song
\cite{CS3}. In \S5 we present some open questions and possible extensions of our results.

\section{Variational formulas}

Unless otherwise explicitly mentioned, we assume throughout this section that 
$\alpha = 1$. 
We briefly recall the connection  between our eigenvalue problem
(\ref{cauchyproblem})  and the mixed Steklov problem discussed in \cite{BK}. Let $D$
be a bounded Lipschitz domain (see \cite{BK} for the precise definition of Lipschitz domain). 
  For
$f
\in L^{1}(\Rd)$ we set
$$P_{t}f(x) = \int_{\Rd} p(t,x,y) f(y) \, dy$$
    where $p(t,x,y)$ is given by (\ref{transition}).
For  $f\in L^{2}(D)$ we extend it to  all of  $\Rd$ by putting $f(x) 
= 0$ for $x
\in
D^{c}$.  Since  $D$ is bounded we see  that such functions are also in
$ L^{1}(\Rd)$.   Thus
$P_{t}f(x)$ is well defined for $f \in L^{2}(D)$ by our bound on $p(t,x,y)$ and
  in particular it is well
defined   for any eigenfunction
$\varphi_n$ of our eigenvalue
problem (\ref{cauchyproblem})  extended to be zero outside of $D$.   
For any $n \in \N$, $x \in \Rd$ and $t > 0$ we put
\begin{equation}
\label{un}
u_{n}(x,t) = P_{t}\varphi_{n}(x) \quad \text{and} \quad u_{n}(x,0) = 
\varphi_{n}(x).
\end{equation}

This defines a function in $$H=\{(x, t):x\in
\Rd, \, t\geq 0\}.$$  
For  bounded Lipschitz domains,  
$\varphi_{n}$ 
is continuous on all of  $\Rd$ (see \cite{BK}, inequality (3.2)), so that $u_{n}$ is 
continuous on all of  $H$.
We will denote by  $H_{+}$ the interior of the set $H$.
  That is, $H_{+}=
\{(x,t): x \in\Rd, t >0\}.$
Let
$$
\Delta = \sum_{i = 1}^{d} \frac{\partial^{2}}{\partial
x_{i}^{2}} + \frac{\partial^{2}}{\partial t^{2}}
$$
denote the Laplace
operator in
$H_{+}$.

We have (Theorem 1.1 \cite{BK})
\begin{eqnarray}
\Delta u_{n}(x,t) & = & 0;
\quad \quad \quad \quad \, \,  \quad (x,t) \in H_{+},
\label{fsteklov1}\\
\frac{\partial u_{n}}{\partial t} (x,0)  & = & - \lambda_{n} u_{n}(x,0);
\quad  \quad x \in  D,
\label{fsteklov2}\\
u_{n}(x,0)  & = & 0;
\quad \quad \quad \quad \quad \, \, \, \quad x \in D^{c}.
\label{fsteklov3}
\end{eqnarray}

The problem (\ref{fsteklov1})--(\ref{fsteklov3}) is called a mixed Steklov problem. The
functions $u_n$ are
called Steklov eigenfunctions. On bounded domains this problem has been 
extensively studied  (see for example,  \cite{HP}, \cite{DS}).
The transformation of our eigenvalue problem (\ref{cauchyproblem}) 
for the Cauchy process to
(\ref{fsteklov1})--(\ref{fsteklov3}) enables us to 
derive a variational formula for
$\lambda_{n}$. This was done in \cite{BK}, Theorem 3.8. 

In this paper we will prove variational formulas for eigenvalue
 gaps $\lambda_n - \lambda_1$
 and for $\lambda_* - \lambda_1$. 
 For $D \subset \Rd$ we set 
$$
H_D = H_+ \cup \{(x,0) \in H: \, x \in D\}.
$$
For $\eps > 0$ we set
$$
H_{\eps} = \{(x,t) \in H: \, t > \eps\}.
$$
By  $\nabla$ we denote the ``full" gradient in $H$.  That is,
$$
\nabla = \left(\frac{\partial}{\partial x_{1}}, \ldots, \frac{\partial}{\partial 
x_{d}}, \frac{\partial}{\partial t} \right).
$$
For brevity, $D_1, \ldots, D_d$ will denote $\frac{\partial}{\partial x_1}, \ldots, \frac{\partial}{\partial x_d}$
and $D_{d+1}$ will denote $\frac{\partial}{\partial t}$. Similarly, $D_1^2, \ldots, D_d^2$ will denote
$\frac{\partial^2}{\partial x_1^2}, \ldots, \frac{\partial^2}{\partial x_d^2}$ and $D_{d+1}^2$ will denote
$\frac{\partial^2}{\partial t^2}$. Coordinate axes in H will be denoted by $0x_1, \ldots, 0x_d, 0x_{d + 1}$ and 
$0x_{d+1}$ denotes the $0t$ axis.

\begin{defn}
\label{piecewiseR1}
We say that a function $f: \R \to \R$ is {\it piecewise $C^1$ on $\R$} 
if the following conditions (i) and (ii) are satisfied:
\begin{itemize}
\item[(i)] There exist a set $A \subset \R$ consisting of at most finitely many points ($A$ may be empty) such
that for any $x \in \R \setminus A$ the derivative $f'(x)$ exists, is finite and continuous at $x$.

\item[(ii)] $f'$ is bounded on $\R \setminus A$.
\end{itemize}
\end{defn}

If we assume that $f:\R \to \R$ is piecewise $C^1$ on $\R$ and $f$
 is continuous on $\R$ then $f$ has the following basic property. For any $a,b \in \R$ we have
$$
\int_a^b f'(x) \, dx = f(b) - f(a).
$$

We shall need the definition of the class of $C^1$ functions on 
 $H_{\eps}$.

\begin{defn}
Let $\eps > 0$ and $f: H_{\eps} \to \R$. We say
 that $f$ is { \it piecewise $C^{1}$ on $H_{\eps}$} if the 
following conditions (i) and (ii) are satisfied for each $i
= 1, \ldots, d, d + 1$. 
\begin{itemize}
\item[(i)] For any line $l \subset H_{\eps}$ parallel to $0x_i$ when $i=1, \ldots, d$ or half-line $l \subset
H_{\eps}$ parallel to $0x_i$ when $i = d + 1$ there exists a subset $A(l,i) \subset l$ ($A(l,i)$ depends on $l$ and
$i$) consisting of at most finitely many points ($A(l,i)$ may be empty) such that for any $(x,t) \in l \setminus
A(l,i)$ the derivative $D_i f(x,t)$ exists, is finite and is continuous at $(x,t)$ as a function on $l$.

\item[(ii)] There exists a constant $c(\eps,i)$ such that for any $(x,t) \in H_{\eps}$ which 
 does not belong to any $A(l,i)$ we have $|D_i f(x,t)| \le c(\eps,i)$.
\end{itemize}
\end{defn}

In the variational formulas for $\lambda_n - \lambda_1$ the   functions $u_n/u_1$ will play a crucial role. 
 We know that $\varphi_1 > 0$ on $D$ so $u_1 > 0$ on $H_D$,
which implies that $u_n/u_1$ is well defined on $H_D$. Since for any $n \in \NP$,  $u_n$ is 
continuous on $H_D$, 
$u_n/u_1$ is also continuous on $H_D$. Intrinsic ultracontractivity for the semigroup $\{P_t^D\}_{t \ge 0}$
proved in \cite{K} implies that that for any $n \in \NP$ there exists a constant $c(D,n)$ such that for
any for any  $x \in D$ we have $\varphi_n(x) \le c(D,n) \varphi_1(x)$. It follows from this that $u_n/u_1$ is
bounded on $H_D$. We also have
\begin{equation}
\label{quotient}
D_i\left(\frac{u_n}{u_1}\right)(x,t) = 
\frac{(D_i u_n(x,t)) u_1(x,t) - (D_i u_1(x,t)) u_n(x,t)}{u_1^2(x,t)}
\end{equation} 
and
\begin{equation}
\label{u1formula}
u_1(x,t) = \int_D \frac{c_d t}{(t^2 + |x - y|^2)^{(d + 1)/2}} \varphi_1(y) \, dy.
\end{equation}
Fix $\eps > 0$. Note that there exists a constant $c(D,\eps)$ such that for any $(x,t) \in H_{\eps}$ and $y \in D$ we have $t^2 + |x - y|^2 \le c(D,\eps) (t^2 + |x|^2)$. 
It follows that there is  a constant $c(D,\eps)$ such that for any $(x,t) \in H_{\eps}$ we have $u_1(x,t) \ge
c(D,\eps) t (t^2 + |x|^2)^{-(d + 1)/2}$. Lemma 3.3(e) in \cite{BK} states that there exists a constant $c(D,n,\eps)$
such that for any $n \in \NP$ and $(x,t) \in H_{\eps}$ we have $|\nabla u_n(x,t)| \le c(D,n,\eps) (t^2 + |x|^2)^{-(d
+ 1)/2}$. Therefore, we see from (\ref{quotient}) that for any $i = 1, \ldots, d + 1$ and $n \in \NP$, $n \ge 2$
the derivative $D_i (u_n/u_1)$ is bounded on $H_{\eps}$. In fact, there exists a constant $c = c(D,n,\eps)$ such that
$\nabla(u_n/u_1)(x,t) \le c/t$ for any $(x,t) \in H_{\eps}$.

We will now introduce the classes of functions $\mathcal{G}(D)$ and $\mathcal{G}_n(D)$
 which we shall use in the variational characterization of $\lambda_n - \lambda_1$. (Note that the set
$\mathcal{G}(D)$ is a linear space.)

\begin{defn}
\label{classG}
Let $D \subset \Rd$ be a bounded Lipschitz domain. We define $\mathcal{G}(D)$ to
 be the set of all functions $u:H_D \to \R$ satisfying the following conditions:
\begin{itemize}
\item[(i)] $u$ is continuous and bounded on $H_D$.

\item[(ii)] For any $\eps > 0$ $u$ is piecewise $C^1$ on $H_{\eps}$.

\item[(iii)]
$$
\int_{H} |\nabla u(x,t)|^2 u_1^2(x,t) dx \, dt < \infty.
$$
\end{itemize}
\end{defn}

When $D \subset \Rd$ is fixed and $u: H_D \to \R$,  we simply set  $\tilde{u}(x) = u(x,0)$, $x \in D$ and 
$$
||\tilde{u}||_2 = \left(\int_D \tilde{u}^2(x) \, dx \right)^{1/2}.
$$

\begin{defn}
\label{classGn}
Let  $D \subset \Rd$ be a bounded Lipschitz domain.  For $n\geq 2$, set 
$$
\mathcal{G}_n(D) = \{u \in \mathcal{G}(D): \tilde{u}\varphi_1 \perp \varphi_1, \ldots, \varphi_{n - 1}; \, ||\tilde{u}\varphi_1||_2 = 1 \}.
$$
\end{defn}
We will often simply write  $\mathcal{G}(D)$  for $\mathcal{G}$ 
and $\mathcal{G}_n(D)$ for $\mathcal{G}_n$
when there is no danger of confusion.

\begin{thm}
\label{variationalf}
Let $D \subset \Rd$ be a bounded Lipschitz domain. Then for any  $n \ge 2$ we have
$$
\lambda_{n} - \lambda_{1} = \inf_{u \in \mathcal{G}_n} \int_{H} |\nabla 
u(x,t)|^{2} u_1^2(x,t) \, dx \, dt.
$$
Moreover, the function $u_n/u_1 \in \mathcal{G}_n$ and the infimum is achieved on this function. 
That is,  
$$
\lambda_{n} - \lambda_{1} =  \int_{H} 
\left|\nabla\left(\frac{u_n}{u_1}\right)(x,t)\right|^{2} u_1^2(x,t) \, dx \, dt.
$$
\end{thm}

\begin{defn}
\label{classG*}
Let $D \subset \Rd$ be a connected bounded Lipschitz domain which
 is symmetric relative to the $x_1$-axis. We set 
$$
\mathcal{G}_*(D) = \{u \in \mathcal{G}(D): \text{ $\tilde{u}$ is
 antisymmetric relative to $x_1$-axis and  $||\tilde{u}\varphi_1||_2 = 1$} \}.
$$
\end{defn}

As above,  we will often write  $\mathcal{G}_*(D)$ for 
$\mathcal{G}_*$. Put $u_*(x,t) = P_t \varphi_*(x)$, 
$(x,t) \in H_+$, $u_*(x,0) = \varphi_*(x)$, $x \in \Rd$
as in formula (\ref{un}).

\begin{thm}
\label{variationalf*}
Let $D \subset \Rd$ be a connected bounded Lipschitz domain
 which is symmetric relative to the $x_1$-axis. We have
$$
\lambda_{*} - \lambda_{1} = \inf_{u \in \mathcal{G}_*} \int_{H} |\nabla 
u(x,t)|^{2} u_1^2(x,t) \, dx \, dt.
$$
Moreover,  the function $u_*/u_1 \in \mathcal{G}_*$ and the 
infimum is achieved on this function.  That is 
$$
\lambda_{*} - \lambda_{1} =  \int_{H} 
\left|\nabla\left(\frac{u_*}{u_1}\right)(x,t)\right|^{2} u_1^2(x,t) \, dx \, dt.
$$
\end{thm}

The proofs of these results 
 will be very similar to the proofs of the variational
formulas for $\lambda_n$ and $\lambda_*$ proved in \cite{BK} 
(see the proofs of Propositions 3.4, 3.5, 3.6, 3.7,
Theorem 3.8 and Proposition 4.8 in \cite{BK}). As in \cite{BK}, 
we first need  some auxiliary propositions.

\begin{prop}
\label{Heps}
Let $D \subset \Rd$ be a bounded Lipschitz domain and assume 
that $u: H_D \to \R$ satisfies conditions (i) and (ii) from Definition \ref{classG}.
 Then for $\eps > 0$ and $n \ge 2$,  
\begin{eqnarray}
&&
\nonumber
\int_{H_{\eps}} \nabla u(x,t) \nabla \left(\frac{u_n}{u_1}\right)(x,t) u_1^2(x,t) \, dx \, dt \\
&&
\label{firstf}
= - \int_{\Rd} u(x,\eps) u_1^2(x,\eps) \frac{\partial}{\partial t} \left(\frac{u_n}{u_1}\right)(x,\eps) \, dx.
\end{eqnarray}
Both integrals are absolutely convergent.
\end{prop}

\begin{proof}
First note that if $f: \R \to \R$ is 
piecewise $C^1$ on $\R$, $g: \R \to \R$ 
 is $C^2$ on $\R$ and  $h: \R \to \R$  is $C^1$ on $\R$, then a simple integration by
parts gives
\begin{equation}
\label{parts1}
\int_a^b f' g' h = [f g' h]_a^b - \int_a^b f g'' h - \int_a^b f g' h',
\end{equation}
for any $a,b \in \R$, $a < b$.
 To prove (\ref{firstf}) we need a multidimensional version of
(\ref{parts1}).  For this we need  some more notation. For any $\eps > 0$, $a > \eps$ let 
$$
\Omega = \Omega(a,\eps) = \underbrace{(-a,a) \times \ldots \times
  (-a,a)}_{d \, \, \text{times}} \times (\eps,a). 
$$ 
Of course, $\Omega \subset H_+$. Let $f: H_+ \to \R$, $g: H_+ \to \R$ and $h: H_+ \to \R$.
Assume that for any $\eps > 0$ $f$ is piecewise $C^1$ on $H_{\eps}$, $g$ is $C^2$ on $H_+$
and $h$ is $C^1$ on $H_+$. Then (\ref{parts1}) implies that for any $\eps > 0$ and any $a >
\eps$ we have
\begin{equation}
\label{parts2}
\int_{\Omega} (\nabla f) (\nabla g) h = \int_{\partial \Omega} f (D_{\nu} g) h - \int_{\Omega} f (\Delta g) h - \int_{\Omega} f (\nabla g) (\nabla h),
\end{equation}
where $D_{\nu}$ is the outer normal derivative on $\partial \Omega$.

The identity (\ref{parts2}) is a well known  version of the Green formula, see
\cite{CH}, page 280, formula 5.  But here, because of a very simple shape of $\Omega$ 
this formula  follows directly from  (\ref{parts1}).

Let us fix $\eps > 0$, $a > \eps$ and apply (\ref{parts2}) to $f = u$, $g = u_n/u_1$, $h =
u_1^2$.  We have 
\begin{eqnarray}
&& 
\label{parts3}
\int_{\Omega} (\nabla u) \left(\nabla \left(\frac{u_n}{u_1}\right)\right) u_1^2 = 
\int_{\partial \Omega} u \left(D_{\nu} \left(\frac{u_n}{u_1}\right)\right) u_1^2 
\\
&&
\nonumber
- \int_{\Omega} u \left( \Delta \left(\frac{u_n}{u_1}\right) \right) u_1^2 - 
\int_{\Omega} u \left( \nabla \left(\frac{u_n}{u_1}\right) \right) (\nabla (u_1^2)) = \text{I} - \text{II} - \text{III}.
\end{eqnarray}
We first  calculate the integrals II and III. Recall that for any $i = 1, \ldots, d, d + 1$ we
have
$$
D_i \left(\frac{u_n}{u_1}\right) = \frac{(D_i u_n) u_1 - (D_i u_1) u_n}{u_1^2}.
$$
Simple calculations gives 
$$
D_i^2 \left(\frac{u_n}{u_1}\right) = \frac{(D_i^2 u_n) u_1 - (D_i^2 u_1) u_n 
- 2 (D_i u_1) (D_i u_n)}{u_1^2} + \frac{2 (D_i u_1)^2 u_n}{u_1^3}.
$$
It follows that 
\begin{eqnarray}
&&
\nonumber
\text{II} = 
\int_{\Omega} u u_1^2 \sum_{i = 1}^{d + 1} D_i^2 \left(\frac{u_n}{u_1}\right) 
\\
&&
\label{fII}
= \int_{\Omega} u u_1 \sum_{i = 1}^{d + 1} D_i^2 u_n
-  \int_{\Omega} u u_n \sum_{i = 1}^{d + 1} D_i^2 u_1
- 2 \int_{\Omega} u  \sum_{i = 1}^{d + 1} (D_i u_1) (D_i u_n)
\\
&&
\nonumber
+ 2 \int_{\Omega} u \frac{u_n}{u_1} \sum_{i = 1}^{d + 1} (D_i u_1)^2.
\end{eqnarray}
Since  the functions $u_n$  are all harmonic in $H_+$, $\Delta u_n = 0$ and it follows 
that the first two integrals in  (\ref{fII}) are zero.

Similarly, 
\begin{eqnarray*}
&&
\text{III} = 
\int_{\Omega} u  \sum_{i = 1}^{d + 1} \left( D_i \left(\frac{u_n}{u_1}\right) \right)
D_i(u_1^2) 
\\
&&
=  2 \int_{\Omega} u  \sum_{i = 1}^{d + 1} (D_i u_1) (D_i u_n)
- 2 \int_{\Omega} u \frac{u_n}{u_1} \sum_{i = 1}^{d + 1} (D_i u_1)^2.
\end{eqnarray*}
Comparing the expressions for II and III we obtain that $\text{II} + \text{III} = 0$.
By (\ref{parts3}) we get
\begin{equation}
\label{secondf}
\int_{\Omega} (\nabla u) \left(\nabla \left(\frac{u_n}{u_1}\right)\right) u_1^2 = 
\int_{\partial \Omega} u u_1^2 \left(D_{\nu} \left(\frac{u_n}{u_1}\right)\right). 
\end{equation}

Next we  estimate $u_1^2(x,t)$. For $(x,t) \in \overline{H_{\eps}}$ (the closure of $H_{\eps}$)
we have
\begin{equation}
\label{estu1}
u_1(x,t) = \int_D \frac{c_d t}{(t^2 + |x - y|^2)^{(d + 1)/2}} \varphi_1(y) \, dy 
\le c(D,\eps) (t^2 + |x|^2)^{-d/2}.
\end{equation}
Hence  $u_1^2(x,t) \le c(D,\eps) (t^2 + |x|^2)^{-d}$. Note also that $u$ satisfies condition (i)
from Definition \ref{classG} so $u$ is bounded on $\overline{H_{\eps}}$. By the
remarks before
Definition \ref{classG}, $\nabla(u_n/u_1)$ is bounded on $\overline{H_{\eps}}$ so
that $D_{\nu}(u_n/u_1)$ is bounded on $\partial \Omega = \partial(\Omega(a,\eps))$,
independently on $a$.

The boundary of $\Omega$  consists of $2(d + 1)$ faces. We denote by $(\partial \Omega)_1$
the face which is a subset of $\partial H_{\eps}$. For any $(x,t) \in \partial \Omega
\setminus (\partial
\Omega)_1$ we have $|x|^2 + t^2 \ge a^2$ so for such $(x,t)$ we have $u_1^2(x,t) \le
c_1(D,\eps) a^{-2 d}$. The measure of $\partial \Omega$ is bounded by $c(d) a^d$. It follows
that
$$
\left| \int_{\partial \Omega \setminus (\partial \Omega)_1} u u_1^2 D_{\nu} \left(\frac{u_n}{u_1}\right) \right| \le c(D,\eps) a^{-d},
$$ 
so when $\eps > 0$ is fixed and $a \to \infty$ this integral tends to $0$.
Note that for $(x,t) \in (\partial \Omega)_1$ we have $D_{\nu} = -
D_{d + 1} = - \frac{\partial}{\partial t}$. It follows that 

$$
\lim_{a \to \infty} \int_{\partial \Omega} u u_1^2 D_{\nu} \left(\frac{u_n}{u_1}\right) 
= \lim_{a \to \infty} \int_{(\partial \Omega)_1} u u_1^2 D_{\nu} \left(\frac{u_n}{u_1}\right)
= - \int_{\partial H_{\eps}} u u_1^2  \frac{\partial}{\partial t}\left(\frac{u_n}{u_1}\right).
$$
The last integral is absolutely convergent by (\ref{estu1}). When $\eps > 0$ is fixed and $a \to \infty$ the set $\Omega$ tends to $H_{\eps}$. Therefore the left hand side of (\ref{secondf}) tends to 
$$
\int_{H_{\eps}} (\nabla u) \left(\nabla \left(\frac{u_n}{u_1}\right)\right) u_1^2,
$$
when $a \to \infty$. When $d \ge 2$ this integral is absolutely convergent by (\ref{estu1}) and by the fact that $\nabla u$ and $\nabla (u_n/u_1)$ are bounded on $H_{\eps}$. When $d = 1$ the last integral is absolutely convergent by (\ref{estu1}), the fact that $\nabla u$ is bounded on $H_{\eps}$ and the fact that $\nabla(u_n/u_1)(x,t) \le c/t$ for $c = c(D,n,\eps)$ and any $(x,t) \in H_{\eps}$.
\end{proof}

\begin{prop}
\label{limiteps}
Let $D \subset \Rd$ be a bounded Lipschitz domain and assume that $u: H_D \to \R$ satisfies conditions (i) and (ii)
 from Definition \ref{classG}. Then for  $n \ge 2$ we have 
$$
\lim_{\eps \to 0^+} \int_{\Rd} u(x,\eps) u_1^2(x,\eps) \frac{\partial}{\partial t} \left(\frac{u_n}{u_1}\right)(x,\eps) \, dx
= - (\lambda_n - \lambda_1) \int_{D} \varphi_n(x) \varphi_1(x) u(x,0) \, dx.
$$
\end{prop}

\begin{proof}
Let $r_n$ be defined as in Proposition 3.1 in \cite{BK}. By Proposition 3.2 (iii) in \cite{BK}
we get 
\begin{eqnarray*}
&&
u_1^2(x,\eps) \frac{\partial}{\partial t} \left(\frac{u_n}{u_1}\right)(x,\eps)
= \frac{\partial u_n}{\partial t}(x,\eps) u_1(x,\eps) 
- \frac{\partial u_1}{\partial t}(x,\eps) u_n(x,\eps)
\\
&&
= (- \lambda_n u_n(x,\eps) + P_{\eps}r_n(x)) u_1(x,\eps)
- (- \lambda_1 u_1(x,\eps) + P_{\eps}r_1(x)) u_n(x,\eps).
\end{eqnarray*}
Since $u$ is bounded we obtain 
$$
\left| \int_{Rd} u(x,\eps) P_{\eps}r_n(x) u_1(x,\eps) \, dx \right|
\le ||u||_{\infty} \int_{\Rd} |P_{\eps}r_n(x)| u_1(x,\eps) dx.
$$
The last integral tends to $0$ as $\eps$ tends to $0^+$ by Proposition 3.5 (formula (3.14)) in \cite{BK}. Exactly in the same way
$$
\left| \int_{\Rd} u(x,\eps) P_{\eps}r_1(x) u_n(x,\eps) \, dx \right|
$$
tends to $0$ as $\eps$ tends to $0^+$.

The only thing which remains is to verify hat 
\begin{equation}
\label{unu1}
\lim_{\eps \to 0^+} \int_{\Rd} u_n(x,\eps) u_1(x,\eps) u(x,\eps) \, dx 
= \int_{D} \varphi_n(x) \varphi_1(x) u(x,0) \, dx.
\end{equation}
Note that $u$ is bounded and $\lim_{\eps \to 0^+} u_n(x,\eps) = 
\varphi_n(x)$, $x \in \Rd$ (recall that $\varphi_n(x) = 0$ for $x \in D^c$). By definition of
$u_n$, for any $x \in \Rd$ and $\eps \in (0,1)$ we have
$$
|u_n(x,\eps)| = |P_{\eps} \varphi_n(x)| \le c(D) ||\varphi_n||_{\infty} (1 + \delta_D(x))^{-d - 1},
$$
where $\delta_D(x) = \dist(x,\partial D)$. Now (\ref{unu1})
 follows by the bounded convergence theorem. 
\end{proof}

\begin{prop}
\label{Hunu1}
Let $D \subset \Rd$ be a bounded Lipschitz domain. Then for $n \ge 2$ we have
$$
\int_{H}  \left| \nabla \left(\frac{u_n}{u_1}\right)(x,t)\right|^2 u_1^2(x,t) \, dx \, dt 
= \lambda_n - \lambda_1.
$$
In particular, we conclude that $u_n/u_1$ satisfies condition (iii) of Definition \ref{classG} and hence $u_n/u_1 \in \mathcal{G}$.
\end{prop}

\begin{proof}
Since $u_n/u_1$ satisfies conditions (i) and (ii) from Definition \ref{classG} we can apply Propositions \ref{Heps} and \ref{limiteps}. This gives
\begin{eqnarray*}
&&
\int_{H}  \left| \nabla \left(\frac{u_n}{u_1}\right)(x,t)\right|^2 u_1^2(x,t) \, dx \, dt
= \lim_{\eps \to 0^+} 
\int_{H_{\eps}}  \left| \nabla \left(\frac{u_n}{u_1}\right)(x,t)\right|^2 u_1^2(x,t) \, dx \, dt
\\
&&
= - \lim_{\eps \to 0^+} \int_{\Rd} \frac{u_n(x,\eps)}{u_1(x,\eps)} u_1^2(x,\eps) \frac{\partial}{\partial t} \left(\frac{u_n}{u_1}\right)(x,\eps) \, dx
\\
&&
=  (\lambda_n - \lambda_1) \int_{D} \varphi_n(x) \varphi_1(x) 
\frac{u_n(x,0)}{u_1(x,0)} \, dx = \lambda_n - \lambda_1.
\end{eqnarray*}
\end{proof}

\begin{prop}
\label{uG}
Let $D \subset \Rd$ be a bounded Lipschitz domain and $u \in \mathcal{G}$. Then for 
 $n \ge 2$ 
$$
\int_{H} \nabla u(x,t) \nabla \left(\frac{u_n}{u_1}\right)(x,t) u_1^2(x,t) \, dx \, dt 
= (\lambda_n - \lambda_1) \int_{D} \varphi_n(x) \varphi_1(x) u(x,0) \, dx. 
$$
Both integrals are absolutely convergent.
\end{prop}

\begin{proof}
Since $u$ and $u_n/u_1$ satisfy condition (iii) of Definition \ref{classG} we have
\begin{eqnarray*}
&&
\lim_{\eps \to 0^+} \int_{H_{\eps}} \nabla u(x,t) \nabla \left(\frac{u_n}{u_1}\right)(x,t) u_1^2(x,t) \, dx \, dt \\
&&
= \int_{H} \nabla u(x,t) \nabla \left(\frac{u_n}{u_1}\right)(x,t) u_1^2(x,t) \, dx \, dt
\end{eqnarray*} 
and the integral on the right hand side is absolutely convergent.  The proposition follows
from Propositions \ref{Heps} and \ref{limiteps}.
\end{proof}

To simplify notation set
$$
Q(u,v) = \int_{H} \nabla u(x,t) \nabla v(x,t) u_1^2(x,t) \, dx \, dt.
$$
Note that for any $u,v \in \mathcal{G}$ the expression $Q(u,v)$ is well defined and finite.

\begin{proof}[Proof of Theorem \ref{variationalf}]
We must show that $\lambda_n - \lambda_1 = \inf_{u \in \mathcal{G}_n} Q(u,u)$. 
Of course, $u_n/u_1 \in \mathcal{G}_n$ and by Proposition \ref{Hunu1},
$$
\inf_{u \in \mathcal{G}_n} Q(u,u) \le Q(u_n/u_1, \, u_n/u_1) = \lambda_n - \lambda_1.
$$

It remains to show that $$\inf_{u \in \mathcal{G}_n} Q(u,u) \ge \lambda_n - \lambda_1.$$
Fix $u \in \mathcal{G}_n$. For  $(x,t) \in H_D$ set
$$
v_{k}(x,t) = (\sum_{m = 1}^{k} c_m u_m(x,t))/u_1(x,t),
$$
where $c_m = \int_{D} \tilde{u}(x) \varphi_1(x) \varphi_m(x) \, dx$.
Since $u \in \mathcal{G}_n$, $n \ge 2$ we know that $\tilde{u} \perp \varphi_1^2$ so $c_1 = 0$. It follows that $v_k \in \mathcal{G}$ because $u_m/u_1 \in \mathcal{G}$ ($m \ge 2$) and $\mathcal{G}$ is a linear space. We have
\begin{equation}
\label{Q(u,u)}
Q(u,u) = Q(v_k,v_k) + Q(u - v_k, u- v_k) + 2 Q(u - v_k, v_k)
\end{equation}
and
\begin{equation}
\label{Quvv}
Q(u - v_k, v_k) = \sum_{m = 1}^{k} c_m Q(u,u_m/u_1) - 
\sum_{m = 1}^{k} c_m Q(v_k,u_m/u_1).
\end{equation}
By Proposition \ref{uG} the right hand side equals 
$$
\sum_{m = 1}^{k} c_m (\lambda_m - \lambda_1) 
\left( \int_{D} \varphi_m(x) \varphi_1(x) u(x,0) \, dx  
- \int_{D} \varphi_m(x) \varphi_1(x) v_k(x,0) \, dx \right).
$$
But $\int_{D} \varphi_m(x) \varphi_1(x) u(x,0) \, dx = c_m$ and for $m = 1, \ldots , k$
$$
\int_{D} \varphi_m(x) \varphi_1(x) v_k(x,0) \, dx = 
\sum_{l = 1}^{k} \int_{D} \varphi_m(x) \varphi_1(x) c_l \varphi_l(x)/\varphi_1(x) \, dx
= c_m.
$$ 
Thus  the expression in (\ref{Quvv}) must be zero. We have also shown that $Q(v_k,v_k) =
\sum_{m = 1}^{k} c_m^2 (\lambda_m - \lambda_1)$. Since $u \in \mathcal{G}_n$ we have
$||\tilde{u} \varphi_1||_2 = 1$, $\tilde{u} \varphi_1 \perp \varphi_1, \ldots, \varphi_{n - 1}$ so
$c_1 = \ldots = c_{n - 1} = 0$. Hence $\sum_{m = n}^{\infty} c_m^2 = 1$. Therefore for $k \ge
n$ we get by (\ref{Q(u,u)})
$$
Q(u,u) \ge Q(v_k,v_k) = \sum_{m = n}^{k} c_m^2 (\lambda_m - \lambda_1)
\ge (\lambda_n - \lambda_1) \sum_{m = n}^{k} c_m^2.
$$
Since $k \ge n$ is arbitrary, we conclude that $Q(u,u) \ge \lambda_n - \lambda_1$.
\end{proof}

\begin{proof}[Proof of Theorem \ref{variationalf*}]
We must show that $\lambda_* - \lambda_1 = \inf_{u \in \mathcal{G}_*} Q(u,u)$. 
Assume that $\lambda_*$ has multiplicity $m \ge 1$ and that it is one of the
eigenvalues
$\lambda_k = \ldots = \lambda_{k + m -1}$, for some $k \ge 2$.  
We may  assume that $u_*
= u_k$.  Note also that $u_*/u_1 \in \mathcal{G}_*$. By Proposition \ref{Hunu1} we get 
$$
\lambda_* - \lambda_1 = \lambda_k - \lambda_1 = Q(u_k/u_1, u_k/u_1)
= Q(u_*/u_1, u_*/u_1) \ge \inf_{u \in \mathcal{G}_*} Q(u,u).
$$
It remains to show that 
$$ \inf_{u \in \mathcal{G}_*} Q(u,u) \ge \lambda_* - \lambda_1.$$

We have 
$\lambda_* - \lambda_1 = \lambda_k - \lambda_1 = 
\inf_{u \in \mathcal{G}_k} Q(u,u)$, where $\mathcal{G}_k = \{u \in \mathcal{G}:
\tilde{u}\varphi_1 \perp \varphi_1, \ldots, \varphi_{k - 1}; \, ||\tilde{u}\varphi_1||_2 = 1 \}$.
By the proof of Proposition 4.8 in \cite{BK}, $\varphi_1, \ldots, \varphi_{k - 1}$ are all
symmetric relative to $x_1$--axis. It follows that $\mathcal{G}_* \subset \mathcal{G}_k$ and
hence, 
$$
\lambda_* - \lambda_1 = \inf_{u \in \mathcal{G}_k} Q(u,u) \le \inf_{u \in \mathcal{G}_*} Q(u,u).
$$ 
\end{proof}

We end this section with   two lemmas which allow us to replace the
Steklov eigenfunction
$u_1(x,t)$ in the variational formula by the simpler expression 
$e^{- \lambda_1 t} \varphi_1(x)$.

\begin{lem}
\label{r1}
For any $x \in D$ and $t > 0$ we have
$$
u_1(x,t) \ge e^{- \lambda_{1} t} \varphi_{1}(x).
$$
\end{lem}
\begin{proof}
By Proposition 3.2 (iii) in \cite{BK} we have
$$
\frac{\partial
u_{1}}{\partial t} (x,t) = -\lambda_{1} u_{1}(x,t)+
P_{t}r_{1}(x),
$$
for $ x \in \Rd, \, t >0$. Moreover we have $r_{1}(x) \ge 0$ for all $x \in \Rd$ 
This follows from Proposition 3.1 \cite{BK} and the fact that $\varphi_1 \ge 0$ on $\Rd$. Put $f(x,t) = e^{ \lambda_{1} t} 
u_{1}(x,t) $, $(x,t) \in H$. We have
$$
\frac{\partial f}{\partial t} (x,t) =  e^{ \lambda_{1} t} P_{t}r_{1}(x) \ge 0,
$$
for $ x \in \Rd, \, t >0$. Hence, for each fixed $x \in D$ the function $f(x,t)$ 
is nondecreasing as a function of $t$. Therefore for $x \in D$ we have  $f(x,t) 
\ge f(x,0) = \varphi_{1}(x)$.
\end{proof}

The following lemma is  an immediate conclusion of Theorem \ref{variationalf*} and Lemma
\ref{r1}. 
\begin{lem}
\label{D01}
We have
$$
\lambda_{*} - \lambda_{1} \ge 
  \int_{0}^{\infty} 
\int_{D} \left| \nabla \left( \frac{u_{*}}{u_{1}} \right)(x,t) \right|^{2} 
\varphi_{1}^{2}(x) e^{-2 \lambda_1 t} \, dx \, dt.
$$
\end{lem}

\section{Weighted Poincar\'e inequalities}

 Let us recall that the positive function $g$ defined on the interval $(-l, l)$ is log--concave if 
the function $\log(g)$ is concave in $(-l, l)$.  That is, for all $x, y\in (-l, l)$ and $0\leq \lambda\leq
1$, 
$$
\log(g(\lambda x+(1-\lambda) y)\geq \lambda \log(g(x)) +(1-\lambda)\log(g(y))
$$
or  equivalently, 
$$
g(\lambda x+(1-\lambda) y)\geq \log(g(x))^{\lambda}\log(g(y))^{1-\lambda}.
$$
If  $g$ is a positive function defined on a convex domain $D\subset \R^d$, then $g$ is said to be
log--concave on $D$ if it is log--concave on every segment contained in $D$.  The celebrated theorem
of Brascamp and Lieb \cite{BrLi} asserts that in the case of Brownian motion, $\varphi_1$ is
log--concave if
$D$ is convex.  In fact, their result is more general than that and it is one of this more general versions 
that we shall use below.  We  state it here in the form that we need. 
Let us recall that in the introduction we have defined (see  (\ref{gaussian}))  
$$
p_t^{(2)}(x)={\frac{1}{(4\pi t)^{d/2}}}e^{-\frac{|x|^2}{4t}}.
$$
This is just  the
Gaussian density in $\R^d$. This is the density for Brownian motion running at twice the usual speed. By $B_t$ we denote the standard Brownian motion in $\Rd$. That is, in our notation we have $P_x(B_{2t} \in A) = \int_{A} p_t^{(2)}(x - y) \, dy$, $x \in \Rd$, $t > 0$, $A \subset \Rd$. 

\begin{prop}\label{BrLi} (Brascamp--Lieb \cite {BrLi}) Let $D \subset \R^d$ be a bounded convex domain
and for
$n \in \N$, let
$t_1, t_2, \dots, t_n$ be  real numbers in $(0, \infty)$. For $x\in D$  define the function 
\begin{equation} G_n(x; \, t_1, \dots, t_n)= \int_D\cdots\int_D\, \prod_{i=1}^n
p_{t_i}^{(2)}(x_{i-1}-x_i)\,dx_1\ldots dx_n,
\end{equation}
where $x_0=x$.  As a function of $x$,   $G_n(x; \, t_1, t_2, \cdots, t_n)$ is log--concave in $D$. 
\end{prop}
Note that
$$
G_n(x; \, t_1, \dots, t_n)=P_x\{B_{2 t_1}\in D,
B_{2(t_1+t_2)}\in D, \dots, B_{2(t_1+t_2+\dots +t_n)}\in D\}.
$$
 
Our desired Poincar\'e inequality will follow from this proposition,  subordination and inequalities
already known for log--concave functions.  First, we recall the latter. 
\begin{prop}
\label{PWS1}
(Payne--Weinberger \cite{PW}, Smits \cite{Sm}). Let $l >0$, $g: (-l,l) \to \R$ be positive and
 log-concave. Let $f: (-l,l) \to \R$ be piecewise $C^1$ and satisfying 
$$ 
\int_{-l}^{l} f(x) g(x) \, dx = 0.
$$
Then 
$$
\int_{-l}^{l} (f'(x))^{2} g(x) \, dx \ge
\frac{\pi^{2}}{4 l^{2}} \int_{-l}^{l} f^{2}(x) g(x) \, dx.
$$
\end{prop}

As an easy consequence of this proposition we obtain the following result.

\begin{cor}
\label{PWS1a}
Let $l >0$, $g: (-l,l) \to \R$ be positive, log-concave, and satisfying  $g(-x) = g(x)$, 
$x \in (-l,l)$.  That is $g$ is symmetric. Let $f: (-l,l) \to \R$ be piecewise $C^1$ and satisfying
$f(-x) = - f(x)$, $x \in (-l,l)$.  That is $f$ is antisymmetric. Then
$$
\int_{-l}^{l} (f'(x))^{2} g(x) \, dx \ge
\frac{\pi^{2}}{4 l^{2}} \int_{-l}^{l} f^{2}(x) g(x) \, dx. 
$$
\end{cor}

From now on we assume that $D$ satisfies the assumptions of 
 Theorem \ref{main} , $L = \sup\{x_1: x=(x_1,\ldots,x_d) \in D\}$.
As an easy conclusion of the above corollary  we get the following 
proposition.

\begin{prop}\label{smits} (Smits\cite{Sm}) 
\label{PWS2}
Let  $g: D \to \R$ be positive, log-concave, and satisfying  $g(\widehat{x}) = g(x)$,
 $x \in D$.  That is, $g$ is symmetric relative to $x_1$--axis. Let $f: D \to \R$, $f \in
C^{\infty}(D)$ and satisfy $f(\widehat{x}) = - f(x)$, $x \in D$. That is, $f$ is antisymmetric
relative to $x_1$--axis. Then 
\begin{equation}
\label{PWS2f} \int_{D} \left| \frac{\partial f}{\partial x_1} (x) \right|^{2} g(x) \, dx \ge 
\frac{\pi^{2}}{4 L^{2}} 
 \int_{D} f^{2}(x) g(x) \, dx. 
\end{equation}
\end{prop}

These type of inequalities are commonly known as
 Poincar\'e inequalities (see Payne-Weinberger \cite{PW}) .

Although we are not able to prove that the first eigenfunction $\varphi_{1}$ 
for the Cauchy process for the domain $D$ is log-concave,  we will be  able  to 
show that the assertion of the previous proposition holds for $g = \varphi_{1}^{2}$ using 
Proposition \ref{BrLi} and subordination. That is, we have 
 
\begin{thm}
\label{PWS3}
Let $f: D \to \R$, $f \in C^{\infty}(D)$ and satisfying  $f(\widehat{x}) = - f(x)$,  
$x \in D$.  That is, $f$ is antisymmetric relative to $x_1$--axis. Then 
\begin{equation}
\label{PWS3f}
\int_{D} \left| \frac{\partial f}{\partial x_1} (x) \right|^{2} \varphi_1^2(x) \, dx \ge 
\frac{\pi^{2}}{4 L^{2}} 
 \int_{D} f^{2}(x) \varphi_1^2(x) \, dx
\end{equation}
where $\varphi_1$ is the first eigenfunction for the symmetric stable process of index
$0<\alpha<2$. 
\end{thm}

\begin{proof}
Let us recall that for  $ 0< \alpha
<2$ the symmetric stable process $X_t$ in $\R^d$ has the representation
\begin{equation}X_t = B_{2\sigma_t},
\end{equation}
where $\sigma_t$ is a  stable subordinator of index $\alpha /2$
independent of $B_t$ (see \cite{BG1}). Thus
\begin{equation}\label{subordination}
p_t^{(\alpha)}(x-y)= \int_0^{\infty} p_s^{(2)}(x-y)
g_{\alpha/2}(t,s)ds,
\end{equation} where $ g_{\alpha/2}(t,s)$ is the
transition density of $\sigma_t$.

Let $x \in D$, $n \in \N$, $0<t_1<t_2 <\dots, t_n$ and set $x_0=x$ and
$t_0=0$.  Using the Markov property for the stable process $X_t$, the subordination
formula (\ref{subordination}), Fubini's theorem, in this order, we obtain,
\begin{eqnarray}
&& F_n(x; \,  t_1, \dots, t_n)=P_x\{X_{t_1}\in D, \dots, X_{t_n}\in
D\}\label{3.6}\\ && = \int_D\cdots\int_D\, \prod_{i=1}^n
p^{(\alpha)}_{t_i-t_{i-1}}(x_{i-1}-x_{i})\,dx_1\ldots
dx_n\nonumber\\ && =\int_0^{\infty}\ldots\int_0^{\infty}\left(
\int_D\cdots\int_D\, \prod_{i=1}^n
p_{s_i}^{(2)}(x_{i-1}-x_{i})\,dx_1\ldots dx_n\right)\\ &&\times
\prod_{i=1}^{n}g_{\alpha/2}(t_i-t_{i-1},s_i)\,ds_1\ldots ds_n\\
&& = \int_0^{\infty}\ldots\int_0^{\infty}G_n(x; \, s_1, \dots, s_n)
\prod_{i=1}^{n}g_{\alpha/2}(t_i-t_{i-1},s_i)\,ds_1\ldots ds_n\label{3.9},
\end{eqnarray}
where $G_n$ is defined as in Proposition \ref{BrLi}.

Let us note that the product of log-concave functions is log-concave.
Using this, Proposition \ref{BrLi} and Proposition \ref{smits}, for each sequence of positive numbers $s_1,
s_2, \dots, s_n$ and $\tilde{s_1}, \tilde{s_n}, \dots, \tilde{s_m}$, $n,m \in \N$  we have (with $f$ as in the statement of the
theorem),
 
\begin{eqnarray*}
&&\int_{D} \left| \frac{\partial f}{\partial x_1} (x) \right|^{2} G_n(x; \, s_1, \dots s_n)G_m(x; \, \tilde
s_1,
\dots\tilde  s_m)
\, dx\\
&& \ge 
\frac{\pi^{2}}{4 L^{2}} 
 \int_{D} f^{2}(x) G_n(x;  \, s_1, \dots s_n)G_m(x; \tilde s_1,
\dots\tilde  s_m)\, dx.
\end{eqnarray*}
Integrating this inequality with respect to $s_1\dots s_n$ and $\tilde s_1, \dots \tilde s_m$ we obtain by 
(\ref{3.6})--(\ref{3.9}),

\begin{eqnarray}\label{stable}
&&\int_{D} \left| \frac{\partial f}{\partial x_1} (x) \right|^{2} F_n(x;\,   t_1, \dots t_n)F_m(x; \, \tilde t_1,
\dots\tilde  t_m)
\, dx\\
&& \ge 
\frac{\pi^{2}}{4 L^{2}} 
 \int_{D} f^{2}(x) F_n(x;\, t_1, \dots t_n)F_m(x; \tilde t_1,
\dots\tilde  t_m)\, dx,\nonumber
\end{eqnarray}
for $0<t_1<t_2<\cdots <t_n$ and $0<\tilde t_1 <\tilde t_2<
\dots<\tilde  t_m$.

Now, let $\tau_D=\inf\{t \ge 0: X_t\notin D\}$. Since $D$ is bounded
and has a Lipschitz boundary Lemma 6 from \cite{Bo} gives that $P_x(X(\tau_D) \in \partial D) = 0$, for any $x \in D$. Using this and the right
continuity of the sample paths we obtain that for any $x \in D$
\begin{eqnarray}\label{approx}
P_x\{ \,\tau_D >t\,\} &=& P_z\{\,  X_{s} \in D, \, \forall \,\,\,   0\leq s
\leq t\,\}\\
&=&\lim_{n \to \infty} P_x \{\, X_{\frac{it}{n}} \in
D, i= 1,\ldots,n\,\}\nonumber\\
&=&\lim_{n \to \infty}
F_n\left(x;\, {\frac{t}{n}} , \, \frac{2t}{n}, \,  \dots \frac{(n-1)t}{n}, \,   t\right) . \nonumber
\end{eqnarray}
Fix $t > 0$, let $n,m \in \N$, $t_i = it/n$, $i = 1, \ldots,n$ and $\tilde t_i =it/m$, $i =1,
\ldots,m$. Letting $n$ and $m$ go to $\infty$,  it follows from (\ref{stable}) and
(\ref{approx}) that

\begin{equation}\label{exittime}
\int_{D} \left| \frac{\partial f}{\partial x_1} (x) \right|^{2}\left( P_x\{ \,\tau_D >t\,\} \right)^2
\, dx \ge 
\frac{\pi^{2}}{4 L^{2}} 
 \int_{D} f^{2}(x)\left(P_x\{ \,\tau_D >t\,\}\right)^2\, dx 
\end{equation}
for all $t>0$.

From the {\it ``intrinsic ultracontractive"} properties of the semigroup
for stable processes in general bounded domains (see \cite{CS1}, \cite{CS2}, \cite{K}), it follows that for any
symmetric stable process 

\begin{equation}\label{egenlimit}
\lim_{t \to \infty}e^{\lambda_1 t}P_x\{ \,\tau_D >t\,\}=\varphi_1(x)
\end{equation}
and this convergence is uniform for $x\in D$.  The inequality (\ref{PWS3f}) follows from (\ref{exittime})
and (\ref{egenlimit}) and the theorem is proved.

\end{proof}

We  call the inequality (\ref{PWS3f}) a ``weighted Poincar\'e--type  inequality for stable processes."
It is interesting to note that the 
eigenfunction $\varphi_1$ in  (\ref{PWS3f}) can
 be replaced by various other similarly
 generated functions from $P_x\{ \,\tau_D >t\,\}$.  For example, we may replace $\varphi_1$ by
$E_x(\tau_D^p)$  or by $\left(E_x \tau_D\right)^p$,  for any $0<p<\infty$.  In addition, the
 theorem holds for any process obtained from Brownian motion by subordination such as the 
relativistic process studied  in \cite{R}.

\section{Proof of Theorems \ref{main} and \ref{upperestimate}}

Unless otherwise explicitly mentioned, we assume throughout this section that 
$\alpha = 1$ and that $D$ satisfies the assumptions of Theorem \ref{main}.  
We shall now apply the results of the previous section and our variational characterization for
$\lambda_*-\lambda_1$ to prove Theorem \ref{main}.

As an immediate conclusion of
Theorem \ref{PWS3} we get the following proposition.
\begin{prop}
\label{PWS4}
Let $u : D \times (0,\infty) \to \R$ be such that for any $t \in (0,\infty)$ the
 function $u(\cdot,t) \in C^{\infty}(D)$. Assume also that $u(\widehat{x},t) = -u(x,t)$ for any
$x \in D$ and $t \in (0,\infty)$.  Then for any $t \in (0,\infty)$ we have
$$
 \int_{D} \left| \frac{\partial u}{\partial x_1}(x,t)\right|^{2} \varphi_{1}^{2}(x) \, dx 
\ge  \frac{\pi^{2}}{4 L^{2}} 
 \int_{D} u^{2}(x,t) \varphi_{1}^{2}(x) \, dx.
$$
\end{prop}
Recall that $\nabla$ is the
``full" gradient in
$H$, that is,  
$$
\nabla = \left(\frac{\partial}{\partial x_{1}}, \ldots, \frac{\partial}{\partial 
x_{d}}, \frac{\partial}{\partial t} \right).
$$
Observe that the function   $u(x,t) = u_*(x,t)/u_1(x,t)$ satisfies the
 assumptions of Proposition \ref{PWS4}. Therefore
\begin{eqnarray}
\label{gradient}
&& \int_{0}^{\infty} \int_{D} \left| \nabla \left( \frac{u_*}{u_1} 
\right) (x,t) \right|^{2} \varphi_{1}^{2}(x) e^{- 2 \lambda_1 t} dx \, dt \\
\nonumber
&& \ge \min \left( \frac{\pi^{2}}{4 L^{2}}, 1 \right) 
\int_{0}^{\infty} \int_{D} \left( \frac{u_*^{2}(x,t)}{u_1^{2}(x,t)} + 
\left| \frac{\partial}{\partial t}\left(\frac{u_*}{u_1}\right)(x,t)\right|^{2}
 \right) \varphi_{1}^{2}(x) e^{- 2 \lambda_1 t} dx \, dt.
\end{eqnarray}

\begin{lem}
\label{derivative}
Let $f : [0,\infty) \to \R$ be a bounded continuous function such that its first derivative
 $f'$ exists and is bounded on $[0,\infty)$. Then for any $c > 0$ we have
\begin{equation}
\label{c+1}
I(f) = \int_{0}^{\infty} (f^{2}(t) + (f'(t))^{2}) e^{-c t} \, dt \ge \frac{f^{2}(0)}{c + 1}.
\end{equation}
\end{lem}
\begin{proof}
We have
$$
I(f) \ge \int_{0}^{\infty} (-2 f(t) f'(t)) e^{-c t} \, dt 
= f^2(0) - c \int_{0}^{\infty} f^{2}(t) e^{-c t} \, dt.
$$
It follows that
$$
(c + 1) I(f) \ge c \int_{0}^{\infty} f^{2}(t) e^{-c t} \, dt + I(f) \ge f^2(0).
$$
\end{proof}
We do not know whether the inequality (\ref{c+1}) is optimal. 
Note only that if we put $f(t) \equiv f(0)$, $t \ge 0$,  then $I(f) = f^2(0)/c$.

\begin{proof}[Proof of Theorem \ref{main}]
By Lemma \ref{D01}, (\ref{gradient}) and Lemma \ref{derivative},   we obtain  
\begin{eqnarray*}
&& \lambda_{*} - \lambda_{1} \ge 
 \int_{0}^{\infty} 
\int_{D} \left|\nabla \left(\frac{u_{*}(x,t)}{u_{1}(x,t)} \right) \right|^{2} 
\varphi_{1}^{2}(x) e^{- 2 \lambda_{1} t} \, dx \, dt\\
&& \ge  \min \left( \frac{\pi^{2}}{4 L^{2}},1 \right) 
 \int_{D}  \int_{0}^{\infty} 
\left( \frac{u_*^{2}(x,t)}{u_1^{2}(x,t)} + 
\left| \frac{\partial}{\partial t}\left(\frac{u_*}{u_1}\right)(x,t)\right|^{2} \right) 
e^{- 2 \lambda_{1} t} \, dt
 \varphi_{1}^{2}(x) \, dx. \\
&& \ge \frac{1}{2 \lambda_1 + 1} \min \left( \frac{\pi^{2}}{4 L^{2}},1 \right) 
\int_{D}  \frac{u_{*}^{2}(x,0)}{u_{1}^{2}(x,0)}  \varphi_{1}^{2}(x) \, dx \\
&&
= \frac{1}{2 \lambda_1 + 1}
\min \left( \frac{\pi^{2}}{4 L^{2}},1 \right)
 \int_{D} \varphi_{*}^{2}(x) \, dx = \frac{1}{2 \lambda_1 + 1}
\min \left( \frac{\pi^{2}}{4 L^{2}},1 \right), 
\end{eqnarray*}
using the fact that 
$$\int_{D} \varphi_{*}^{2}(x) \, dx =1.$$
Rewriting this we find that 
\begin{equation}
\label{final}
\lambda_* - \lambda_1 \ge  
 \min \left( \frac{\pi^2}{4  (2 \lambda_1 + 1) L^2}, \frac{1}{2 \lambda_1 +1} \right).
\end{equation}

Since the inradius of D is equal to $1$ we have $\lambda_1  \le \lambda_1(B(0,1))$. 
By Corollary 2.2 \cite{BK} we have $\lambda_1(B(0,1)) \le C(d)$, where
$$
C(d) = \frac{\pi d (d +2)}{4 (d + 1)}.
$$
It follows that $\lambda_1 \le C(d)$.  This and 
 (\ref{final}) conclude the proof of Theorem \ref{main}. 
\end{proof}

\begin{proof}[Proof of Theorem \ref{upperestimate}] We first recall the following deep result
of Ashbaugh and Benguria,  the so called ``Payne--P\'olya--Weinberger
conjecture" proved in  \cite{AB2} .   For any bounded connected domain $D\subset \Rd$, we denote by
$\mu_2(D)$ and $\mu_1(D)$ the second and first eigenvalues of the
Dirichlet Laplacian in $D$, respectively. (Of course, $\mu_2(D)$ and
$\mu_1(D)$  are the second and first eigenvalues for the semigroup of
Brownian motion killed upon exiting $D$.) Let $B$ be any ball in $\Rd$.  
The Payne--P\'olya--Weinberger conjecture proved
in \cite{AB2} asserts that 

\begin{equation}\label{PPW}
\frac{\mu_2(D)}{\mu_1(D)}\leq
 \frac{\mu_2(B)}{\mu_1(B)}=\frac{j_{{d/2}, 1}^2}{j_{{d/2-1},1}^2}
\end{equation}
Furthermore,  equality holds if and only if $D$ is a ball (we will  not use this fact here).    
To avoid confusion let us also denote by $\lambda_1(D)$ and $\lambda_2(D)$ the first and second eigenvalues
for the semigroup of the Cauchy  process killed upon exiting $D$.  It follows by the upper bound in \cite{BK} and the
lower bound in \cite{CS3} that for $i=1, 2$ and for convex domains $D$,  
\begin{equation}\label{gapupper}
\frac{1}{2} \sqrt{\mu_i(D)}\leq \lambda_i(D)\leq \sqrt{\mu_i(D)}.
\end{equation}
From this, 
\begin{equation}
\lambda_2(D)-\lambda_1(D)\leq \sqrt{\mu_2(D)}-{\frac{1}{2}}\sqrt{\mu_1(D)}.
\end{equation}
However, by (\ref{PPW}), 
$$
\frac{\sqrt{\mu_2(D)}}{{\frac{1}{2}}\sqrt{\mu_1(D)}}-1\leq
\frac{\sqrt{\mu_2(B)}}{{\frac{1}{2}}\sqrt{\mu_1(B)}}-1,
$$
where here we choose $B$ to be the largest ball contained in the domain $D$. 
  This inequality can be written as 
$$
{\frac{\sqrt{\mu_2(D)}-{\frac{1}{2}}\sqrt{\mu_1(D)}}{{\frac{1}{2}}\sqrt{\mu_1(D)}}}\leq 
{\frac{\sqrt{\mu_2(B)}-{\frac{1}{2}}\sqrt{\mu_1(B)}}{{\frac{1}{2}}\sqrt{\mu_1(B)}}}
$$
which leads to 
\begin{eqnarray}
\sqrt{\mu_2(D)}-{\frac{1}{2}}\sqrt{\mu_1(D)}&\leq&
\left(\sqrt{\mu_2(B)}-{\frac{1}{2}}\sqrt{\mu_1(B)}\right)
\left(\frac{\sqrt{\mu_1(D)}}{\sqrt{\mu_1(B)}}\right)\\
&\leq & \sqrt{\mu_2(B)}-{\frac{1}{2}}\sqrt{\mu_1(B)},\nonumber
\end{eqnarray}
where we used the fact that 
$\sqrt{\mu_1(D)}\leq \sqrt{\mu_1(B)}$, by domain monotonicity of the first eigenvalue. 
By scaling, 
$\mu_2(B)=\frac{j_{{d/2}, 1}^2}{r_D^2}$ and $\mu_1(B)=\frac{j_{{d/2-1}, 1}^2}{r_D^2}$,
which proves  the desired inequality.  
\end{proof}

Let $\lambda_i(D)$ be the eigenvalues for the semigroup of the symmetric $\alpha$-stable process killed on
exiting a bounded convex domain $D$. Using the more general inequality 
\begin{equation}
\frac{1}{2} \left(\mu_i(D)\right)^{\alpha/2}\leq \lambda_i(D)\leq \left(\mu_i(D)\right)^{\alpha/2},
\end{equation}
valid for any $0 < \alpha < 2$, see \cite{DeMe},  \cite{CS3} and the argument
above we have the following generalization of Theorem \ref{upperestimate}.
\begin{thm}
Let $D\subset \Rd$ be a bounded convex domain of inradius $r_D$.  Let $\lambda_1$ and
$\lambda_2$ be the  first and second eigenvalues for the semigroup of the symmetric 
$\alpha$-stable process $0 < \alpha < 2$ killed upon
exiting $D$.  Then 
\begin{equation} 
\lambda_2 - \lambda_1 \leq
{\frac{j_{{d/2}, 1}^{\alpha}-{\frac{1}{2}}j_{{d/2-1}, 1}^{\alpha}}{r_D^{\alpha}}}.
\end{equation}
\end{thm}

In the case of Brownian motion  the above argument gives  that for any bounded domain $D$ of
inradius $r_D$, 
\begin{equation} 
\mu_2(D) - \mu_1(D) \leq \mu_2(B) - \mu_1(B)=
{\frac{j_{{d/2}, 1}^{2}-j_{{d/2-1}, 1}^{2}}{r_D^{2}}}
\end{equation} 
with equality if and only if $D$ is a ball.    We believe the following conjecture
should be true. 
\begin{conj} Let $D\subset \Rd$ be a bounded domain and let $\lambda_2(D)$ and $\lambda_1(D)$ be the
second and first eigenvalues for the semigroup of the symmetric $\alpha$-stable process $0 < \alpha < 2$ killed upon
exiting $D$.  Then 
\begin{itemize}
\item[(i)] (The $\alpha$--stable version of the Payne--P\'olya--Weinberger Conjecture): 
\begin{equation*} 
\frac{\lambda_2(D) }{\lambda_1(D)} \leq
\frac{\lambda_2(B) }{\lambda_1(B)} 
\end{equation*}
with equality if and only if $D$ is a ball.  In particular, 
\begin{equation*} 
\lambda_2(D) - \lambda_1(D) \leq \lambda_2(B) - \lambda_1(B)
\end{equation*} 
with equality if and only if $D$ is a ball.
\item [(ii)] If $D$ has inradius $r_D$, then 
\begin{equation*} 
\lambda_2(D) - \lambda_1(D) \leq
{\frac{j_{{d/2}, 1}^{\alpha}-j_{{d/2-1}, 1}^{\alpha}}{r_D^{\alpha}}}.
\end{equation*}
\end{itemize}
\end{conj}
We refer the reader to \cite{BLM} and  \cite{Me} where many of the classical
isoperimetric--type inequalities which hold for Brownian motion are shown to also hold for
symmetric stable processes. 

As for a conjecture concerning a sharp lower bound we have the following (see also Remark
\ref{remark5.1} below).

\begin{conj}
Let $D\subset \Rt$ be a bounded convex domain which is symmetric relative to both
coordinate  axes.  Let $R$ be the smallest oriented (sides parallel to the coordinate axes) rectangle
containing $D$.  For any $0<\alpha<2$,
\begin{equation}
\lambda_2(R) - \lambda_1(R) \leq \lambda_2(D) - \lambda_1(D).
\end{equation}
\end{conj}

For  Brownian motion ($\alpha=2$) this is proved in \cite{BKr}, \cite{BM}, \cite{Da}.
\section{Concluding Remarks}

We end this paper with several remarks and questions  which 
naturally arise from our results.

\begin{rem}\label{remark5.1}
  For planar domains $D$  with the symmetry assumptions  of  Theorem
\ref{main} and for  Brownian motion,  it follows from \cite{BKr}, \cite{BM}, \cite{Da} 
that $\lambda_* - \lambda_1
\ge 3\pi^2/ (4 L^2)$,  and for arbitrary convex domains of diameter $d$,
$\lambda_2-\lambda_1 > \pi^2/d^2$, 
\cite{Sm}, \cite{YZ}.  We may ask whether our estimates  for the Cauchy process, 
$\lambda_* -
\lambda_1$  is optimal in terms of the order of  $L$.   Let us
consider the rectangle
$R = [-L,L] \times [-1,1]$
where $L \ge 1$. In a forthcoming paper we will show  that
\begin{equation}
\frac{c\ln(L+1)}{L^{2}} \leq  \lambda_{*} - \lambda_{1} \leq  \frac{C \ln(L+1)}{L^{2}}
\end{equation}
for two absolute positive constants $c,C$.  For this case, the methods in this paper give only
$1/L^{2}$  due to the fact that we integrate the expression in the variational
formula over $D
\times [0, \infty)$ (see Lemma \ref{D01}) and the extra term $\ln(L+1)$  comes from 
integration over all of $H$.

When  $D \subset \Rt$ is a dump--bell shaped  domain (say two disjoint unit balls
joined by a sufficiently thin corridor) 
 which is symmetric according to the $x_1$-axis,  one can show that 
\begin{equation}
\lambda_{*} - \lambda_{1} \le \frac{C}{L^{3}},
\end{equation}
where $C > 0$ does not depend on $D$ and $L >> 1$.
(Since  trivially $\lambda_2
-\lambda_1 \le \lambda_* - \lambda_1$,  the upper bound estimate for
$\lambda_* - \lambda_1$  also gives  the same estimate for the spectral gap
$\lambda_2 - \lambda_1$.)
Thus  the lower bound result of this paper is not true for arbitrary
non-convex domain.  It may also be that we have here a different
situation
 than in the case of Brownian motion case where 
the spectral gap $\lambda_2 - \lambda_1$ tends to zero as 
 the corridor becomes thinner and thinner and the domain becomes two
 disjoint balls.
It is probably the case that  the spectral gap $\lambda_2 - \lambda_1$ (for the Cauchy
process) of this dump-bell tends to the spectral gap of the set 
which consists of two disjoint balls, and the spectral gap for such a
set is strictly positive.

The  existence and
 properties of $\lambda_*$ and $\varphi_*$ (Theorem 4.3 \cite{BK}) were formulated and
proved for connected, bounded and symmetric Lipschitz domains.  In
fact these assumptions were needed only for technical reasons and the
existence
and other basic properties   are true without the assumptions of connectedness
and Lipschitz boundary. 
This leads to the following question. Assume $D\subset \R^2$ has diameter $d_D$, inradius
$r_D$ and is  symmetric relative to the
$x_1$--axis.    What is the best lower bound
 estimate for $\lambda_*
-\lambda_1$ in terms of $d_D$ and $r_D$ (regardless of connectedness or convexity of $D$)? 
Of course, the same question may be asked for the spectral gap $\lambda_2 - 
\lambda_1$. These questions are  non-trivial even in the one-dimensional
case when $D$  consists  of  finite  number of disjoint
intervals. 
\end{rem}

\begin{rem}
It may be possible to apply the techniques used in this paper and in \cite{BK} to study
eigenvalues and eigenfunctions  for other processes.  Of course, the most obvious 
extensions would be to other symmetric stable processes.  It would also 
be of interest  to extend  these results 
to the relativistic process (\cite{CMS}, \cite{R}) with characteristic
function
 $E^0 e^{i \xi X_t} = e^{- t
(\sqrt{m^2 + |\xi|^2} - m)}$, $t > 0$, $\xi \in \Rd$, $m > 0$.
The infinitesimal generator of this
process is the so-called relativistic Hamiltonian $- \sqrt{- \Delta + m^2} + m$. As with  the
Cauchy process one can build a ``relativistic" Steklov problem  of the
type
\begin{eqnarray*}
&&
\Delta u_n(x, t)+ 2 m \frac{\partial u_n}{\partial t}(x,t)  =  0;
\quad \quad    \, \,  \quad (x,t) \in H_{+},
\\ &&
\frac{\partial u_n}{\partial t} (x,0)   =  - \lambda_{n} u_n(x,0);
\quad  \quad \quad \quad \quad \, x \in  D,
\\ &&
u_{n}(x,0)   =  0;
\quad \quad \quad \quad \quad \quad \quad \quad \quad \, \, \, \, 
\quad x \in D^{c}.
\end{eqnarray*}
Using the identity  
$$
e^{2 m t} \left( \Delta + 2 m \frac{\partial}{\partial t} \right)
= e^{2 m t} \left( \sum_{i = 1}^{d} \frac{\partial^2}{\partial x_i^2} \right) 
+ \frac{\partial}{\partial t} \left( e^{2 m t} \frac{\partial}{\partial t} \right) 
$$
one can show that the eigenvalues of the relativistic process are given by the variational 
formula 
$$
\lambda_{n} = \inf_{u \in \tilde{\mathcal{F}}_n} \int_{H} |\nabla 
u(x,t)|^{2} e^{2 m t} \, dx \, dt,
$$
for an appropriately chosen class of functions $\tilde{\mathcal{F}}_n$.
Thus  the eigenvalue problem for the relativistic process 
is similar to that of  the Cauchy process. Nevertheless, extending the results which we now have for the Cauchy
process remains mostly open  (although some results follow from   the recent paper 
 \cite{CS3}, see Example 6.2). 
\end{rem}

\begin{rem}
As mentioned in the introduction, the spectral gap $\lambda_2-\lambda_1$ measures the
rate at which the  Cauchy process conditioned to remain forever in the domain $D$ tends
to equilibrium.  That is, for any $\varepsilon>0$, we define (as in \cite{SaCo}) the time to
equilibrium $T_{\varepsilon}$ by 
\begin{equation}
T_{\varepsilon}=\inf\{t>0: \sup_{x,y \in D} \left|\frac{e^{\lambda_1 t}p_D(t,x,y)}{ \varphi_1(x)
\varphi_1(y)} - 1
 \right| \leq \varepsilon\}.
\end{equation}
It follows from (\ref{spectralgap}) that 
$$
\frac{1}{\lambda_2-\lambda_1}\log{\frac{1}{\varepsilon}}\leq T_{\eps}
\leq
C_1+\frac{1}{\lambda_2-\lambda_1}\log{\frac{1}{\varepsilon}}.
$$

While a probabilistic interpretation of $\lambda_*-\lambda_1$ is not as ``clean" and useful
as the one above, we do have the following.  Recall that $D_+ = \{x \in D: x_1 > 0\}$ and
$D_- =
\{x
\in D: x_1 < 0\}$.  
Then for any $x \in D_+$
\begin{equation}
-(\lambda_*-\lambda_1)=\lim_{t\to\infty}{\frac{1}{t}}
\log\left(\frac{P^{x}(X_t \in D_+ \, , \, \tau_D > t) -
P^{x}(X_t
\in D_-
\, ,
\,
\tau_D
  > t) }{P^{x}(\tau_D > t)}\right). 
\end{equation}
This follows from the proof of Theorem 4.3 in \cite{BK}, the definition
of $\tilde{p}_D(t,x,y)$ (see Lemma 4.5 in \cite{BK}) and the general
theory of semigroups.

\end{rem}

\end{document}